\documentclass [reqno,12pt]{article}

\newtheorem{teo}{Theorem}[section]
\newtheorem{prop}{Proposition}[section]
\newtheorem{lemma}{Lemma} [section]

\newtheorem{cor}{Corollary}[section]
\newtheorem{definiz}{Definition} [section]
\usepackage{amscd}
\usepackage[intlimits]{amsmath} 
\usepackage{amssymb}
\newcommand{\barPG}{P^0_G(M,\om)} 
\newcommand{\PG}{P_G(M,\om)}
\newcommand{\fine} {\hfill\textsc{Q.D.E.}\medskip\\ }
\newcommand{\dimo}[1][] {\noindent\textbf{Proof#1}. }
\newcommand{\eps} { \varepsilon} 
\renewcommand{\phi} {\varphi}

\newcommand{\R} { \mathbb{ R } } 
\newcommand{\C} { \mathbb { C } } 
\newcommand{\Zeta} { \mathbb{Z} } 
 
\newcommand{\I} {\operatorname{\mathrm {i} }} 

\newcommand{\di} {\, d } 
\renewcommand{\det } { \operatorname {det} }
\newcommand{\est} {\raisebox{.34ex} {$\scriptstyle {\bigwedge}$}}
\newcommand{\Aut} {\operatorname {Aut}} 
\newcommand{\Gl}{\operatorname {Gl}}

\newcommand{\vol}{\operatorname{\mathrm {vol}}} 
\newcommand{\Keler} {K\"{a}hler }
 
\newcommand{\barz}{\bar{z}}
\newcommand{\KE} {K\"{a}hler-Einstein }

\newcommand{\debar } {\bar {\partial } }
\newcommand{\Ric}{\operatorname {Ric} } 
\newcommand{\cinf}{C^{\infty}} 
\newcommand{\cinfm} {\cinf (M) } 
 
\newcommand{\OO} {\mathcal{O}} 
\newcommand{\chern}{\operatorname{\mathrm{c}}} 
\newcommand{\om}{\omega} 
\newcommand{\ra} {\rightarrow}
 
\newcommand{\demi} { \frac {1} {2}} 
 
\newcommand{\vacuo} {\emptyset}
\newcommand{\perdef } { : \, = } 
\newcommand{\idd}{\I \partial \debar}

\newcommand{\lamma}{\lambda}
\newcommand{\intv}{\frac{1}{V}\int_M} 
\newcommand{\Fo}{F^0}

\newcommand{\ciclico}[1]{\Zeta_{#1}} 
\newcommand{\PP}{\mathbb{P}}
\newcommand{\barP}{P^0} 

\newcommand{\dio}{\langle}
\newcommand{\rospo}{\rangle}
\newcommand{\MV}{\dio [\om]^n, [M]\rospo}

\begin{document}

\title{Symmetries, Quotients \\ and\\ K\"{a}hler-Einstein metrics}

\author{Claudio Arezzo, Alessandro Ghigi, Gian Pietro Pirola}

\maketitle

\abstract{
\noindent We consider Fano manifolds $M$ that admit a collection of finite automorphism groups $G_1, ..., G_k$, such that the quotients $M/G_i$ are smooth Fano manifolds possessing a \KE metric. Under some numerical and smoothness assumptions on the ramification divisors, we prove that $M$ admits a \KE metric too.
}

\tableofcontents

\section{Introduction}

The aim of this paper is to provide new examples of \KE metrics of
positive scalar curvature.  The existence of such a metric on a Fano
manifold is a subtle problem, due to the presence of obstructions,
that have been discovered during the years, beginning with
Matsushima's theorem in 1957, Futaki invariants in 1982, Tian's
theorem stating that \KE manifolds of positive scalar curvature are
semistable (see \cite[Theorem 8.1]{tian-KE97}), up to Donaldson's
result \cite[Corollary 4]{donaldson-scalar-embedding-I}, which shows
that the existence of \KE metrics (even more generally of a \Keler
constant scalar curvature metric) forces the algebraic underlying
manifolds to be asympotically stable (see also
\cite{arezzo-loi-CMPh-1}).

Existence theorems on the other hand are always very hard.  The only
necessary and sufficient condition, established by Tian, is of a truly
analytic character.  It says that a Fano manifold $M$ admits a \KE
metric, if and only if an integral functional $F$ defined on \Keler
metrics in the class $\chern_1(M)$ is proper (see Theorem
\ref{tian-properness} below).  The equivalence of properness of $F$
with the algebraic stability of the underlying manifold, in an
appropriate sense, would represent the final solution of the problem,
but is still unknown. (This has been suggested by Yau, and made
precise by Tian, who has also proved that properness implies
stability.)  Work in progress by Paul and Tian \cite{double-tian-paul}
indicates a new stability condition as a candidate for the equivalence
with the existence of a \KE metric.

Although by now there is a good deal of examples, the only broad class
of manifolds for which the problem is solved is the one of toric Fano
manifolds, thanks to a recent theorem of Xujia Wang and Xiaohua Zhu
(\cite{wang-zhu}, see also Donaldson's work \cite{donaldson-torico}
for related results for extremal metrics). Otherwise, even for
manifolds that are deceptively simple from the algebro-geometric point
of view, one has often no clue on how to check the properness of $F$,
and finding the metric. The case of Del Pezzo surfaces is quite
eloquent from this point of view, as the reader of
\cite{tian-KE-surfaces-1990} might verify.  Another striking example
of the difficulties on which one suddenly runs, is the hypersurface
case.  Indeed, it is expected that any smooth Fano hypersurface has a
\KE metric, nonetheless the only ones for which this is known are the
ones lying in a suitable small analytic neighbourhood of the Fermat's
hypersurfaces (see \cite [p. 85-87]{tian-libro}).  In fact a standard
implicit function theorem argument shows that the \KE condition is
open in the moduli space in the analytic topology, provided the
automorphism group is finite. This remark can be applied also to some
of the examples discussed below.

In trying to construct explicit examples a good help comes from having
many holomorphic symmetries to work with. This has been crucial for
example in estimating the so called $\alpha$-invariant for some Del
Pezzo surfaces with reductive automorphisms group.  This has been the
heart of the work of Tian-Yau \cite[Proposition
2.2]{tian-yau-surfaces}. 

The aim of this paper is to use in a different way the symmetries of
the manifold to prove existence of \KE metrics, inspired by Tian's work on Fermat hypersurfaces.  In Section
\ref{brutta} we study the behaviour of properness of $F_\om$ (see p. 
\ref{paggina}) in presence of a Galois covering and find conditions
under which the existence of a \KE metric on the base allows one to
prove a version of properness, and thus existence, on the covering
space.   We find algebraic conditions on the covering maps
(Theorems \ref{criterio} and \ref{criterio2}) ensuring that the
desired inequalities hold on the covering space. In
\ref{esempiozzi-lisci} we show how this can be used to prove the
existence of \KE metrics on some classes of Fano manifolds, chosen
from the lists of Del Pezzo manifolds, and Fano threefolds with
$\operatornamewithlimits{Pic}=\Zeta$ (see \cite[p. 
214-215]{fano-varieties}). Our examples include:
\begin{enumerate}
\item[a)] hypersurfaces of the form $ \{ x_0^d + ... +x_{k-1}^d +
  f(x_k , ..., x_{n+1}) = 0 \} \subset \PP^{n+1} $ where $f$ is a
  homogeneous polynomial of degree $d$, and $k > n+2 -d$;
\item[b)] $n$-dimensional intersections of hypersurfaces of the same
  form as above, all of the same degree $d$ and with $k > n+2 -d$;
\item[c)] arbitrary intersections of two (hyper)quadrics;
\item[d)] double covers of $\PP^n$ ramified along a smooth
  hypersurface of degree $2d$ with $\frac{n+1}{2} < d\leq n$;
\item[e) ] double covers of the $n$-dimensional quadric $Q_n \subset
  \PP^{n+1}$ with smooth branching locus cut out by a hypersurface of
  degree $2d$ with $\frac{n}{2} < d < n$. 
\end{enumerate}
(See section \ref{esempiozzi-lisci}.)  Example (a) generalises Tian's
theorem about Fermat's hypersurfaces.  
Examples (a), (b), (d) and (e) give positive-dimensional
\emph{algebraic families} of \KE manifolds.  This becomes even more
striking in example (c) since every element in the moduli of such
manifolds has a \KE metric. A particular case of (c) (the intersection
of two specific quadrics in $\PP^5$) had been previously studied by
Alan Nadel (see \cite[p.  589] {nadel-annals-multiplicator}).

Some interesting questions arise naturally from these results. In the
first place, when a finite group $G$ acts on an algebraic manifold
$M$, the quotient $M/G$ can always be endowed with the structure of a
complex analytic orbifold. We believe that our theorems can be
generalised to cover this case, provided the quotient admits a \KE
\emph{orbifold} metric.  Nevertheless there are few examples of \KE
orbifolds (see e.g. \cite{demailly-kollar-exponent},
\cite{johnson-kollar-fourier}, \cite{boyer-galicki-kollar-Annals}), and it
is probably hard to apply our results to coverings with orbifold base.

From a different
perspective, in light of our results (c)-(f), one could study the
Weil-Petersson geometry of the moduli spaces of these new families, or
one can try to generalise Mabuchi and Mukai's results
(\cite{mabuchi-mukai}) on compactification of moduli spaces. A
situation which seems geometrically appealing is the one of the
intersection of two quadrics (which is in fact Mabuchi-Mukai's case in
dimension 2).  A classical result says the moduli space of the
intersection of two quadrics in $\PP^{2n + 3}$ is isomorphic to the
moduli space of hyperelliptic curves of genus $n$ (see
\cite{avritzer-lange-communications-algebra} and reference therein). 
Therefore this moduli space inherits two Weil-Petersson geometries,
one coming from the \KE metrics on the intersection of quadrics, the
other from Poincar\'e metrics on curves. It would be interesting to
compare them. 

We wish to thank Gang Tian for many helpful conversations and for his
interest in this work., and the referees for useful suggestions.

\section[Existence theorems on covering spaces]{Existence theorems on covering\\ spaces}
\label{brutta}

Let 
$M$ be a compact $n$-dimensional \Keler manifold and $\om$ a smooth
closed (1,1)-form on $M$
such that
\begin{equation*}
  \dio [\om]^n , [M]\rospo = \int_M \om^n >0. 
\end{equation*}
When $\phi\in\cinfm$ put $\om_\phi=\om +\idd \phi$. 
Define \label{pagina-quatto}
the following functionals on $\cinfm$:
\begin{align}
  I_\om(\phi)&=\frac{1}{\MV}\int \phi (\om^n - \om^n_\phi)
  \label{def-di-I}
  \\ J_\om(\phi) 
& =  \int_0^1 \frac{I_\om(s\phi)}{s}\di s
\label{eq:def-di-J}
\\
 \Fo_\om(\phi) & = J_\om(\phi) - \frac{1}{\MV}\int \phi \om^n. 
\end{align}
When no confusion is possible, we will write $V=\MV$. For the reader's
convenience we recall the following equivalent definitions of these
functionals.
\begin{lemma}
  If $M$ and $\om$ are as above, and $\phi \in \cinfm$, then
  \begin{align}
    \label{eq:J-stile-Donaldson}
    J_\om(\phi) & =- \frac{n!}{\MV} \sum_{p=1}^n \frac{1}{(n-p)! 
      (p+1)!} \int_M \phi \om^{n-p} (\idd\phi)^p=\\
    \label{eq:J-stile-Tian}
& =- \frac{1}{\MV} \sum_{k=0}^{n-1} \frac{k+1}{n+1}
  \int_M \phi \idd \phi \wedge \om^k \wedge \om_\phi^{n-k-1}  =\\
\label{eq:J-stile-Tian-per-parti}
&= \frac{1}{\MV} \sum_{k=0}^{n-1} \frac{k+1}{n+1}
  \int_M \I \partial\phi \wedge\debar \phi \wedge \om^k \wedge
  \om_\phi^{n-k-1}  
\intertext{and}
\label{eq:F0-stile-Donaldson}
F_\om^0(\phi) &=  - \frac{n!}{\MV} \sum_{p=0}^n \frac{1}{(n-p)! 
       (p+1)!} \int_M \phi \om^{n-p} (\idd\phi)^p. 
  \end{align}
\end{lemma}
\dimo[ (sketch)] To prove \eqref{eq:J-stile-Donaldson} expand
$\om_{s\phi}^n=(\om + s\idd \phi)^n $ in powers of $s$ and
use the result to compute  $I_\om(s\phi)$ in
\eqref{eq:def-di-J}. 
As for \eqref{eq:J-stile-Tian} compute $I_\om(s\phi)$ in
\eqref{eq:def-di-J} using
the fact that
\begin{gather*}
  \om^n -\om_{s\phi}^n = (\om - \om_{s\phi}) \sum_{q=0}^{n-1}
  \om^{n-q-1}\wedge \om_{s\phi}^q. 
\end{gather*}
Substituting $\om_{s\phi}=s\om_\phi +(1-s)\om$ and expanding 
$\om_{s\phi}^q$ yields
\begin{gather*}
    J_\om(\phi)  =- \frac{1}{\MV} \sum_{p=0}^{n-1} 
C_p \int_M \phi \idd \phi \wedge \om^{n-1-p} \wedge \om_\phi^{p}   
\end{gather*}
where
\begin{gather*}
  C_p = \sum_{q=p}^{n-1} \binom{q }{p} \int_0^1 s^{p+1} (1-s)^{q-p}\di s. 
\end{gather*}
This can be computed using the combinatorial identities
\begin{gather*}
  \int_0^1 s^{p+1}(1-s)^k\di s = \frac{(p+1)! k! }{(p+k+2)!} \\
\sum_{k=0}^{n-p-1} \frac{p+1}{(p+k+1)(p+k+2)} = \frac{n-p}{n+1}
\end{gather*}
and gives the desired result.  To get
\eqref{eq:J-stile-Tian-per-parti} it is enough to integrate by parts,
using that $\om$ is closed and $M$ is \Keler. Finally
\eqref{eq:F0-stile-Donaldson} is an immediate consequence of
\eqref{eq:J-stile-Donaldson}.  \fine 
Formula
\eqref{eq:F0-stile-Donaldson} says that $F^0$ coincides (up to a
constant factor) with the functional called $I$ by other authors. 
Compare with eq. (25) in \cite{donaldson-symmetric} where Donaldson
gives a nice geometric interpretation of $F^0$. 
\begin{lemma}\label{lemma-triciclico}
  Let  $M$ and $\om$ be as above. If $\lambda$ is a positive constant then
  \begin{equation}
    \label{eq:scaling-di-F0}
    F^0_{\lamma\om}(\lamma \phi) = \lamma F^0_\om( \phi). 
  \end{equation}
Let $\om_0$ be a closed (1,1)-form such that $\dio [\om_0], [M]\rospo
>0$. 
Given 
 $\phi_{01}$, $ \phi_{12}\in \cinfm$
put
 $\om_1=\om_0 + \idd \phi_{01}$,
$\phi_{02}=\phi_{01} + \phi_{12}$. Then
\begin{equation}
  \label{eq:triciclo}
  \Fo_{\om_0}(\phi_{02}) = \Fo_{\om_0}(\phi_{01}) +\Fo_{\om_1}(\phi_{12}). 
\end{equation}
\end{lemma}
For the proof see \cite[pp. 60f]{tian-libro}. 

Assume from now on that $M$ is a Fano manifold and $\om$ is a \Keler metric
in the class $2\pi\chern_1(M)$. Then $V=\MV=n!\vol(M)$. Let $f=f(\om)$ be the unique function
on $M$ satisfying
\begin{equation}
  \Ric(\om) = \om + \idd f(\om), \qquad \int_M e^{f(\om) } \om^n = V. 
\end{equation}
Define $A_\om, F_\om : \cinfm \ra \R$ by
\begin{align*}
  A_\om(\phi) &= \log \biggl [ \intv e^{f(\om) - \phi}\om^n \biggr]
\\
  F_\om(\phi )& = \Fo_\om(\phi) - A_\om(\phi). 
\end{align*}
Although these functionals (as well as the ones defined before) are
defined on the whole of $\cinfm$, their interest for \KE metrics lies
in their behaviour on a smaller space, whose definition we now recall. 

 Let $G$ be a
compact group of isometries of  $(M,\om)$. Put
\begin{equation}
  P_G(M, \om) =\{ \phi \in \cinfm : \om_\phi >0, \text{ and $\phi$ is
  $G$-invariant}\}. 
\end{equation}
By $\om_\phi >0$ we mean that $\om_\phi$ is a \Keler metric.  If
$G=\{1\}$ we simply write $P(M,\om)$. 
We say that $F_\om$ is \emph{proper} on $\PG$ if there is a proper
increasing function $\mu : \R \ra \R$, such that the inequality
\label{paggina}
\begin{equation*}
  F_\om(\phi) \geq \mu\bigl (J_\om(\phi)\bigr )
\end{equation*}
holds for any $\phi \in \PG$.  The importance of this notion is mainly
due to the following theorem (see \cite[Theorem 1.6]{tian-KE97} and
\cite[Chapter 7] {tian-libro}). 
\begin{teo}[Tian]\label{tian-properness}
  Let $M$ be a Fano manifold, $G $ a maximal compact subgroup of
  $\Aut(M)$ and $\om$ a $G$-invariant \Keler metric in the class $2\pi\chern_1(M)$. Then $M$ admits a \KE metric if and only if $F_\om$ is
  proper on $P_G(M, \om)$. Moreover, in this case $F_\om$ is bounded
  from below on all $P(M,\om)$. 
\end{teo}
The elements of $\PG$ parametrise metrics only up to a constant,
because $\om_\phi$ does not change by adding a constant to $\phi $,
and the functional $F_\om$ depends on $\phi\in \PG$ only up to a
constant.  Therefore we can normalise the elements of $\PG$ one way or
another.  The following normalisation is useful in this context:
\begin{equation}
  Q_G(M, \om) = \{ \phi \in P_G(M, \om) : A_\om(\phi) = 0\}. 
\end{equation}
For any $\phi\in \PG$, $\phi +A_\om(\phi) \in Q_G(M,\om)$ is the
corresponding normalised potential.

The following proposition gives a sufficient condition for the
existence of \KE metrics on Fano manifolds. 
\begin{prop}
  \label{criterio-del-sup}
 Let $M$ be a Fano manifold,  $\om$ a \Keler metric in the class
  $2\pi\chern_1(M)$ and $G$ a compact group of isometries of $(M,\om)$.  If there are constants $C_1 , C_2 > 0$ such that
  \begin{equation}
    F_\om(\phi) \geq C_1 \sup_M\phi - C_2
  \label{eq:prop1}
  \end{equation}
  for any $\phi \in Q_G(M, \om)$, then $M$ admits a \KE metric. 
\end{prop}
\dimo One exploits the same estimates used in the proof of Theorem
\ref{tian-properness} (compare \cite[Chapter 7]{tian-libro}.)  Indeed,
let $\phi_t$, $t\in [0,T)$ be the curve of potentials obtained by
applying the continuity method:
\begin{equation}
\label{continuita}
  (\om + \idd \phi_t )^n = e^{f -t \phi_t} \om^n. 
\end{equation}
Then it is known that for some constants $C_3 , C_4 > 0$
\begin{gather}
  \Fo_\om(\phi_t) \leq 0
  \label {eq:Fo-negativo}
  \\ F_\om(\phi_t) \leq -A_\om(\phi_t) \leq \frac{1-t}{V} \int_M
  \phi_t \om_t^n\\ 0\leq -\inf_M \phi_t \leq C_3 \biggl ( \intv
  (-\phi_t) \om_t^n +C_4 \biggr ) \\ \intv \phi_t \om_t^n \leq C_4 \\
  F_\om(\phi_t) \leq - A_\om(\phi_t) \leq C_4 (1-t) \leq C_4
\label{eq:stimo-F-continuita}
\end{gather}
(see \cite[p. 72]{tian-libro}).  
Since 
$\phi_t
+ A_\om(\phi_t) \in Q_G(M,\om)$, and $F_\om(\phi_t)$ does not change by
adding a constant to $\phi_t$, an application of 
\eqref{eq:prop1} yields
\begin{equation}
  \label{eq:criterio-primo-su-P}
 \begin{gathered}
   F_\om(\phi_t) =
   F_\om\bigl (\phi_t + A_\om(\phi_t)\bigr)
   \geq \\
   \geq C_1 \sup_M \bigl (\phi_t + A_\om(\phi_t) \bigr ) -C_2 
   =C_1 \sup_M \phi_t +
   C_1 A_\om(\phi_t) -C_2
\end{gathered}
\end{equation}
Therefore
using \eqref{eq:stimo-F-continuita}
\begin{gather*}
C_1 \sup_M \phi_t \leq   F_\om(\phi_t) - C_1 A_\om(\phi_t) + C_2 \leq
C_4 + C_2 +C_1 C_4
\end{gather*}
Hence $\sup_M \phi_t $ is uniformly bounded. But from \eqref {eq:Fo-negativo}
\begin{equation*}
  J_\om(\phi_t) \leq \intv \phi_t \om^n \leq \sup_M \phi_t. 
\end{equation*}
So $J_\om(\phi_t)$ is bounded and this is enough to bound the $C^0$
norm (see \cite[p. 67]{tian-libro}).  Therefore, by Yau's estimates,
one can solve equations \eqref{continuita} up to $t=1$, and $\om +
\idd \phi_1$ is the \KE metric.  \fine
\begin{lemma}
  \label{superlemma-alpha}
  Let $M$ be a Fano manifold, and $\om$ a \Keler metric in the class
  $2\pi\chern_1(M)$.  
Then for any $\beta >0$ there are  constants $C_1,  C_2 >0$ such
that for any $\phi\in Q(M,\om)$
 \begin{equation}
   \label{eq:log-limita-sup}
 \log \biggl [ \intv e^{- (1+\beta)
   \phi} \om^n \biggr ]  \geq C_1 \sup_M \phi - C_2. 
 \end{equation} 
\end{lemma}
\dimo
According to one of the basic results of   Tian's theory of the  $\alpha$-invariant (see
\cite[Prop. 2.1]{tian-certain}),   there are $\alpha
\in (0,1)$ and $C_3 >0$, such that for any $\phi \in \PG$
\begin{equation}
  \intv e^{ - \alpha (\phi - \sup \phi) } \om^n \leq C_3. 
  \label{invariante-alpha}
\end{equation}
Let $p $ be such that
\begin{equation}
  \frac{p-\alpha}{p-1} = 1+\beta. 
\end{equation}
Then
\begin{equation}
  p=1 + \frac{1-\alpha}{\beta}
\end{equation}
so $p\in (1, +\infty)$, because $\alpha < 1$. 
Let 
$\di \mu$ denote the measure $(1/V)\om^n$ on $M$. 
 By definition, if $\phi
\in Q_G(M,\om)$
\begin{equation}
   \int  e^{f-\phi}\di \mu =1,
\label{uguale1}
\end{equation}
so 
$$
e^{-\sup f} \leq \int e^{-\phi} \di \mu. 
$$
Since
\begin{gather*}
  -\phi = \frac{\alpha}{p} (  \sup \phi -\phi) -
 \frac{\alpha}{p} \sup \phi + \biggl ( 1 - \frac{ \alpha}{p} \biggr
 ) (- \phi)
\end{gather*}
\begin{gather*}
e^{-\sup f} \leq e^{-\frac{\alpha}{p} \sup \phi} \int e^{\frac{\alpha}{p}
  ( \sup \phi -\phi) } \cdot e^{\frac{p-\alpha}{p}( - \phi) }\di \mu. 
\end{gather*}
Therefore applying H\"{o}lder inequality with exponent $p$ yields
\begin{equation}
  e^{\frac{\alpha}{p} \sup \phi - \sup f} \leq \biggl [ \int
  e^{\alpha(\sup \phi-\phi) } \di \mu \biggr]^{1/p} \biggl [ \int
  e^{\frac{p'}{p} (p-\alpha) (-\phi)} \di \mu \biggr ]^{1/p'}. 
\end{equation}
Using \eqref{invariante-alpha}  and observing that
\begin{equation*}
  \frac{p'}{p}(p-\alpha) = 1+\beta,
\end{equation*}
we get
\begin{gather*}
  e^{\frac{\alpha}{p} \sup \phi -\sup f} \leq C_3^{1/p} \,
 \biggl [ \intv
  e^{-(1+\beta)\phi}  \om^n \biggr ]^{1/p'}. 
\end{gather*}
Taking logarithms
\begin{gather*}
\frac{p'\alpha}{p} \sup \phi -p'\sup f - \frac{p'}{p} C_3 \leq \log
\biggl [ \intv
  e^{-(1+\beta)\phi}  \om^n \biggr ]
\end{gather*}
that is \eqref{eq:log-limita-sup} with
$$
C_1 = \frac{p'\alpha}{p}= \frac{\alpha \beta}{1-\alpha} > 0, \quad
C_2 = p'\sup f + \frac{C_3}{p-1}. 
$$
\fine
\begin{cor}
\label{criterio-del-beta}
If there are constants $C_1, C_2>0$ and $\beta >0$ such that
 \begin{equation}
   F_\om(\phi) \geq C_1 \log \biggl [ \intv e^{ - (1+\beta)
   \phi} \om^n \biggr ] -C_2
 \end{equation} 
 for any $\phi \in Q_G(M, \om)$, then $M$ admits a \KE metric. 
\end{cor}
This is an immediate consequence of the previous lemma and Proposition
\ref{criterio-del-sup}.

In the proof of the existence theorems below we will need a slight
extension of the integral functionals defined above.  Let $M$ be a
compact complex manifold and $\gamma$ a continuous hermitian form on
$M$. A closed positive current $T$ of bidegree (1,1) is called a
\emph{\Keler current} if for some constant $c>0$ one has $T\geq c
\gamma$ in the sense of currents. The definition does not depend on
the choice of $\gamma$, since $M$ is compact. If $M$ is a Fano
manifold, $G \subset \Aut(M)$ is a compact subgroup, and $\om$ is a
$G$-invariant \Keler form in the class $2\pi \chern_1(M)$, we put
\begin{equation*}
  P^0_G(M,\om) = \{ \psi \in C^0(M) : \om + \idd \psi \text { is a
  \Keler current}\}. 
\end{equation*}
This means that $\psi $ belongs to $P^0_G(M,\om)$ if and only if $\om
+ \idd \psi \geq c \om$ in the sense of currents for some $c>0$. 
\begin{lemma}\label{lemma-estensione}
  The map $ \phi \mapsto (\om + \idd \phi)^n$ can be extended to a
map
\begin{equation*}
  \barPG \longrightarrow \{ \text {positive Borel measures on } M \}. 
\end{equation*}
The extension is continuous with respect to the $C^0$-topology on the
  domain and the weak convergence of measures on the target. 
\end{lemma}  
\dimo This follows from basic results on the complex Monge-Amp\`ere
operator.  Consider a covering
$\{ U_k\}$ of $M$ with
contractible open subsets.  On $U_k$ we have $\om=\idd u_k$ for some
smooth strictly plurisubharmonic function $u_k$.  If $\phi \in \barP$
then $u_k + \phi$ is plurisubharmonic and continuous on $U_k$. 
Although in general currents cannot be multiplied, Bedford and Taylor
showed how to define consistently $\bigl (\idd (u_k + \phi) \bigr )^n$
as a positive measure on $U_k$. Moreover, it follows from the
Chern-Levine-Nirenberg inequality that this measure depends
continuously on $\phi$ (see e.g. \cite[Corollary
2.6]{demailly-Monge-Ampere-operators}).  As $\bigl ( \idd (u_k + \phi)
\bigr)^n= \bigl ( \idd (u_j + \phi) \bigr)^n$ on $U_k\cap U_j$, these
local measures glue together, and the resulting measure on $M$,
denoted by $( \om + \idd \phi)^n$ depends continuously on $\phi$. 
\fine
\begin{prop}\label{prop-estensione}
  The functionals $I_\om$, $J_\om$, $\Fo_\om$ and $F_\om$ can be
  extended to $\barPG$.  The extensions are continuous with respect to
  the $C^0$-topology. 
\end{prop}
\dimo It follows from the previous lemma that we can extend
continuously $I_\om$. Using formula \eqref{eq:def-di-J} we can extend
continuously $J_\om$, and therefore $\Fo_\om$.  $A_\om(\phi)$ can be
clearly extended continuously to $\barPG$.  \fine

In  the proof of the next Theorem we will need the following density result. 
\begin{prop} \label{almost-richberg}
  Any $\psi \in \barPG$ is the $C^0$-limit of a sequence $\phi_n\in
  P_G(M,\om)$. 
\end{prop}
This is  a straightforward 
  application of a result due to Richberg (\cite{richberg}) that we quote in the version given by Demailly (\cite[Lemma 2.15]{demailly-regularization}). 
\begin{lemma}[Richberg]
Let $\psi \in C^0(M)$ be such that $\idd \psi \geq \alpha$ for some
continuous (1,1)-form $\alpha$. Then given any hermitian form $\gamma$
and any $\eps>0$, there is a function $\psi'\in \cinfm$ such that
$\psi \leq \psi' < \psi + \eps $ and $\idd \psi' \geq \alpha -
\eps\gamma$. 
\end{lemma}

The following two lemmata deal specifically with coverings. 
\begin{lemma}\label{lemma-sulle-mappe-finite-corenti}
  If $\pi : M \ra N$ is a finite holomorphic map of compact complex
  manifolds, the direct image via $\pi$ of a \Keler current on $M$ is
  a \Keler current on $N$. 
\end{lemma}
\dimo Let $R\subset M$ and $B\subset N$ denote
ramification and branching locus of $\pi$, and $d$ its degree. Let
$\gamma_M$ and $\gamma_N$ be continuous hermitian forms on $M$ and $N$
respectively. Since $\pi^*\gamma_N$ is continuous and $\gamma_M$ is
positive definite, there is $c_1 > 0$ such that $\gamma_M \geq c_1
\pi^* \gamma_N$. If $T$ is a \Keler current on $M$, by definition
$T\geq c_2 \gamma_M$ for some $c_2 >0$, so that $T\geq c
\pi^*\gamma_N$ with $c=c_1 c_2>0$.  Given a positive form $\eta \in
\est^{n-1,n-1}(N)$ we have
\begin{multline*}
  \langle \pi_*T, \eta \rangle = \langle T , \pi^* \eta \rangle \geq c
  \langle \pi^*\gamma_N, \pi^* \eta \rangle =c \int_M \pi^* ( \gamma_N
  \wedge \eta) =\\ = c\int_{M\setminus R} \pi^* ( \gamma_N \wedge
  \eta) = c \cdot d \int_{N\setminus B} \gamma_N \wedge \eta = c\cdot
  d \int _N \gamma_N \wedge \eta = c\cdot d \langle \gamma_N, \eta
  \rangle
\end{multline*}
so that $ T \geq c \cdot d \gamma_N.$ This proves the lemma.  \fine{}
\begin{lemma}\label{F0-nei-rivestimenti}
  Let $\pi: M \ra N$ be a degree $d$ covering between $n$-dimensional
  \Keler manifolds. Let $\om_N$ be a \Keler metric on $N$, and 
 $\psi \in P^0(N,\om_N)$ a \emph{continuous} potential such that $\pi^*\psi$ be a \emph{smooth}
  function on $M$. Then
  \begin{equation}
    F^0_{\pi^*\om_N} (\pi^*\psi)= F^0_{\om_N}(\psi). 
  \end{equation}
\end{lemma}
\dimo
Put $V_N=\langle N, [\om_N]^n \rangle$. Then $\langle M, [\pi^*
\om_N]^n \rangle = d\cdot V_N$. 
\begin{align}
\begin{split}
  I_{\om_N}(s\psi) & = \frac{1}{V_N} \int_N s\psi \bigl [ \om_N ^n -
  (\om_N + s \idd \psi)^n \bigr ] = \\ & = \frac{1}{dV_N} \int_M
  s\pi^*\psi \bigl [ (\pi^*\om_N) ^n - (\pi^* \om_N +s \idd \pi^*
  \psi)^n \bigr ] =\\ &=  I_{\pi^* \om_N}
  (s\pi^*\psi)
\end{split}
\notag \\
\begin{split}
  J_{\om_N}(\psi) &= \int_0^1 \frac{ I_{\om_N}(s\psi) } {s} \di s =
\int_0^1 \frac{I_{\pi^* \om_N} (s\pi^*\psi)
  }{s} \di s = \\ &=  J_{\pi^* \om_N} (\pi^*
  \psi)
\end{split}
\label{ectoplasma2}
\\
\begin{split}
  \frac{1}{V_N}\int_N \psi \om_N^n & = \frac{1}{dV_N} \int_M \pi^*\psi
  (\pi^*\om_N)^n . 
\end{split}
\label{ectoplasma1}
\end{align}
 Plugging \eqref{ectoplasma2} and \eqref{ectoplasma1} in the definition
of $F^0$ we get finally
\begin{gather}
\Fo_{\om_N} (\psi) 
=\Fo_{\pi^* \om_N} (\pi^*\psi ). 
\label{primo-pezzo}
\end{gather}
Mark that the functionals $I_{\om_N}(\psi), J_{\om_N}(\psi)$ and $ F^0_{\om_N}(\psi)$  are well-defined
because $\om_N$ is a \Keler metric and $\psi \in P^0(N,\om_N)$. 
On the other hand,  $\pi^*\om_N$ degenerates along the ramification,
so it is not a \Keler metric. Nevertheless it is a smooth closed
$(1,1)$-form  and    $\pi^*\psi$
is a smooth function, so   the functionals
$I_{\pi^*\om_N}(\pi^*\psi), J_{\pi^*\om_N}(\pi^*\psi)$ and
$F^0_{\pi^*\om_N}(\pi^*\psi)$ are well-defined too, thanks to the
discussion at p. \pageref{pagina-quatto}. 
\fine
\begin{teo}\label{birazionale-quousque}
  Let $M$ and $N$ be Fano manifolds, $\pi : M \ra N$ a ramified Galois
  covering of degree $d$ with structure group $G$, $\om_N$ a \KE
  metric on $N$ and $\om \in 2\pi\chern_1(M)$ a $G$-invariant \Keler
  metric.  Denote by $R(\pi) $ be the ramification divisor of $\pi$
   (with multiplicities),
  and assume that numerically (i.e. in homology)
\begin{equation*}
  R(\pi) = - \beta K_M
\end{equation*}
for some $\beta \in \mathbb{Q}$. (Since $R(\pi)$ is effective and
$-K_M$ is ample, $\beta > 0$.)  Then there is a constant constant
$ C $ such that for any $\phi \in P_G(M, \om_M)$
\begin{equation}
  F^0_{\om} (\phi) \geq  \frac{1}{1+\beta}
\log \biggl [ \intv e^{- (1+\beta)\phi} \pi^*
  \om^n_N \biggr ] -C. 
\end{equation}
\end{teo}
(Here $V=\MV$.)\\
\dimo The classical Hurwitz formula for the canonical bundle of a
ramified covering, $ \pi^* K_N=K_M - R(\pi)$, yields that
\begin{equation*}
  \pi^*[\om_N]= ( 1 + \beta) [\om]. 
\end{equation*}
Choose a
$G$-invariant $u\in \cinfm$ such that 
\begin{equation}
  \label{eq:def-di-u}
 \pi^*\om_N = (1+\beta) \om +
\idd u. 
\end{equation}
We claim that any $\phi \in P_G(M, \om)$ is of the form
$$
\phi = \frac{u + \pi^* \psi}{1+\beta}
$$
for some
$\psi \in \barP(N,\om_N)$.  Indeed $(1+\beta) \phi -u$ is
$G$-invariant, so $(1+\beta)\phi-u=\pi^*\psi$ for some continuous
function $\psi$, because $N=M/G$ has the quotient topology. 
Since (as currents) 
\begin{equation*}
  \om_N + \idd \psi = \frac{1+\beta}{d} \pi_* ( \om + \idd \phi ),
\end{equation*}
Lemma \ref{lemma-sulle-mappe-finite-corenti} implies that $\om_N +\idd \psi$ is a \Keler
current, i.e. that $\psi \in \barP(N, \om_N)$. We have shown that to
any potential $\phi \in P_G(M, \om)$ corresponds a \emph{continuous}
potential $\psi \in \barP(N, \om_N)$ such that
\begin{equation} 
  \label{definizione-di-psi}
  \pi^* ( \om_N + \idd \psi) = (1+\beta) ( \om + \idd \phi ). 
\end{equation}
Since $N$ is \KE by hypothesis, Tian's Theorem \ref{tian-properness}
implies that there is a constant $C_3$ such that $F_{\om_N} (\eta )
\geq - C_3$ for any $\eta \in P(N, \om_N)$. By Proposition
\ref{prop-estensione} the functional $F_{\om_N}$ can be extended
continuously to $\barP(N, \om_N)$, and by Proposition
\ref{almost-richberg} $P(N, \om_N) $ is dense in $\barP(N, \om_N)$, so
we can conclude that
\begin{equation}
  \label{eq:diseq-su-psi}
 F_{\om_N} (\psi ) \geq - C_3
\end{equation}
for $\psi$ as in \eqref{definizione-di-psi}.  To finish the proof we
need to ``lift'' this inequality from $N$ to $M$. 
From \eqref{eq:def-di-u} and 
\eqref{eq:triciclo} of Lemma \ref{lemma-triciclico}, applied to the
forms $(1+\beta)\om$ and $\pi^*\om_N$, it follows that
\begin{gather*}
\Fo _{(1+\beta)\om} \bigl ( (1+\beta)
  \phi \bigr ) 
=\Fo_{ (1+\beta)\om} (u) + \Fo _{\pi^*\om_N}  (
\pi^*\psi). 
\end{gather*}
Since $u$ does not depend on $\phi$, $\Fo_{ (1+\beta)\om} (u)$ is
a constant. Next  \eqref{eq:scaling-di-F0}
in Lemma \ref{lemma-triciclico} implies that
$$
\Fo _{(1+\beta)\om} \bigl ( (1+\beta)\phi\bigr )= (1+\beta) \Fo_\om (\phi) . 
$$
So
\begin{gather*}
  F^0_\om (\phi) =  \frac{1}{1+\beta} F^0_{\pi^*\om_N}(\pi^*\psi)- C_4. 
\end{gather*}
Using Lemma \ref{F0-nei-rivestimenti} and \eqref{eq:diseq-su-psi} we get
\begin{gather*}
  F_\om^0(\phi) =  \frac{1}{1+\beta} F^0_{\om_N}(\psi) -C_4 \geq \\
\geq \frac{1}{1+\beta}
\log \biggl [
  \frac{1}{V_{N}}\int_N e^{-\psi}\om^n_N \biggr ]
-C_5.
\end{gather*}
Now
\begin{gather*}
\begin{split}
  \frac{1}{V_N}\int_N e^{-\psi}\om_N^n&= \frac{1}{dV_N}\int_M
  e^{-\pi^*\psi}(\pi^*\om_N)^n\\
d V_N &= \dio \pi^* [\om_N^n], [M]\rospo 
=(1+\beta)^n V \\
  \frac{1}{V_N}\int_N e^{-\psi}\om_N^n&=\frac{1}{(1+\beta)^n V}\int
  e^{-(1+\beta)\phi} e^u (\pi^*\om_N)^n \geq \\
&\geq
\frac{e^{\inf u}}{(1+\beta)^n V}\int
  e^{-(1+\beta)\phi}  (\pi^*\om_N)^n.\\
\end{split}
\end{gather*}
Therefore
\begin{gather*}
  F_\om^0(\phi) \geq 
 \frac{1}{1+\beta}
\log \biggl [\intv
  e^{-(1+\beta)\phi}  (\pi^*\om_N)^n \biggr ] -C_6. 
\end{gather*}
\fine
The first criterion for the existence of \KE metrics is the
following
\begin{teo}\label{criterio}
  Let $M$ be an $n$-dimensional Fano manifold. Assume that
  ramified coverings $\pi_i : M \ra M_i$ are given for $i=1, ...,k$,
  satisfying the following assumptions:
   \begin{enumerate}
   \item $M_i$ is a Fano manifold and admits a \KE metric;
   \item  the coverings are Galois, i.e. $M_i=M/G_i$ for some
     finite group $G_i$,
   \item the groups $G_i$ are contained in some \emph{compact}
     subgroup $G \subset \Aut(M)$;
   \item if $R(\pi_1), ..., R(\pi_k)$ are the ramification divisors,
     then
     \begin{equation*}
       \bigcap_{i=1}^k R(\pi_i) = \vacuo;
     \end{equation*}
   \item the divisors $R(\pi_i)$ are all proportional to the
     anticanonical divisor of $M$, i.e. there are some (necessarily
     positive) rational numbers $\beta_i$ such that numerically (i.e. 
     in homology)
     \begin{equation*}
       R(\pi_i) = -\beta_i K_M. 
     \end{equation*}
   \end{enumerate}
   Then $M$ has a \KE metric. 
\end{teo}
\dimo Fix a $G$-invariant \Keler form $\om \in 2\pi \chern_1(M)$ and
\KE metrics $\om_i$ on $M_i$.  As $G_i \subset G$ we have that
\begin{equation}
  P_G(M,\om) \subset \bigcap_{i=1}^k P_{G_i}(M, \om). 
\end{equation}
From Theorem \ref{birazionale-quousque} it follows that for some
constants $C_{1i} \in \R$ we have
\begin{equation}
  F^0_\om(\phi) \geq \frac{1}{1+\beta_i}\log \biggl [ \intv e^{ - (1+\beta_i) \phi}
  (\pi_i ^* \om_i)^n \biggr ] - C_{1i}
\label{eq:trentini}
\end{equation}
for all $ \phi \in P_G(M, \om)$.  Put 
\begin{gather*}
  \begin{gathered}
C_1=\max C_{1i} \\
    \beta=\min \beta_i
  \end{gathered}
\qquad
 \begin{gathered}
p_i=1 + \beta_i
\\
p=\min p_i = 1 +\beta
  \end{gathered}
\qquad 
\begin{gathered}
\psi=e^{-\phi}\\
\di \mu = \frac{1}{V} \om^n
  \end{gathered}
\end{gather*}
and define $\eta_i\in\cinfm$ by
\begin{gather*}
  \pi_i^*\om_i^n=\eta_i\om^n. 
\end{gather*}
Clearly $\beta>0, p>1$ and $\eta_i\geq 0$. Then \eqref{eq:trentini} becomes
\begin{gather*}
  F^0_\om(\phi) +C_{1i} \geq
\frac{1}{p_i} \log \biggl [ \int_M \psi^{p_i} \eta_i \di \mu \biggr]
\end{gather*}
so
\begin{gather*}
  F^0_\om(\phi) +C_1 \geq 
 \log || \psi \eta_i^{1/p_i} ||_{p_i} 
\end{gather*}
where $||\ ||_{s}$ denotes the norm of $L^s(M,\mu)$. 
By construction, $\,p_i\geq p$, for  $i=1,...,k$. If $p_i=p$, then clearly
\begin{gather*}
  || \psi \eta_i^{1/p_i} ||_{p_i} \geq C_{2i} || \psi \eta_i^{1/p} ||_{p} 
\end{gather*}
with $C_{2i}=1$. 
If $p_i>p$, then
\begin{gather*}
  \frac{1}{p}=\frac{1}{p_i} + \frac{1}{q}
\end{gather*}
for some $q>p>1$. By H\"{o}lder inequality
\begin{gather*}
  || \psi \eta_i^{1/p} ||_{p} \leq
|| \psi \eta_i^{1/p_i} ||_{p_i} \cdot ||\eta_i^{1/q}||_q 
\end{gather*}
so
\begin{gather*}
  || \psi \eta_i^{1/p_i} ||_{p_i} \geq C_{2i}
  || \psi \eta_i^{1/p} ||_{p}  
\end{gather*}
with 
$$
C_{2i}=\frac{1}{||\eta_i^{1/q}||_q }>0. 
$$
Actually $||\eta_i^{1/q}||_q =p_i^{n/q}$ so $C_{2i}
 = p_i^{-n/q}$. 
Put $C_2=\min C_{2i}>0$. Then
\begin{gather}
F^0_\om(\phi) +C_1 - \log C_2 \geq \log || \psi \eta_i^{1/p} ||_{p}
=\notag
\\
\notag
=\frac{1}{p}\log\biggl [ \int_M \psi^p \eta_i \di \mu \biggr ]\\
\exp \bigl ( p  F^0_\om(\phi) + C_3 \bigr ) 
\geq \int_M \psi^p \eta_i \di \mu. 
\notag
\end{gather}
Taking the average over $i$
\begin{gather*}
  \exp \bigl ( p   F^0_\om(\phi) + C_3 \bigr ) 
\geq \frac{1}{k}\sum_{i=1}^k\int_M \psi^p  \eta_i \di \mu. 
\end{gather*}
Taking the logarithm we get
\begin{gather}
\notag
F^0_\om(\phi) \geq \frac{1}{p} \log\biggl [\frac{1}{k}\sum_{i=1}^k \int_M \psi^p
\eta_i  \di\mu\biggr ] -C_4  \\
  F^0_\om(\phi) \geq \frac{1}{1+\beta} \log\biggl [ \intv e^{-(1+\beta)\phi}
  \biggl ( \frac{1}{k}\sum_{i=1}^k\eta_i\biggr )\om^n\biggr ] -C_4. 
\label{eq:compare}
\end{gather}
It follows from assumption (4)  that for some
constant $C_5 >0$
\begin{equation}
  \frac{1}{k}\sum_{i=1}^k \eta_i \geq C_5. 
  \label{eq:from-ass-4}
\end{equation}
Therefore
\begin{gather*}
  F^0_\om(\phi) \geq \frac{1}{1+\beta} \log\biggl [ \intv e^{-(1+\beta)\phi}
  \om^n\biggr ] -C_6 . 
\end{gather*}
This holds for any $\phi\in P_G(M\om)$. If $\phi\in Q_G(M,\om)$, then $F_\om(\phi)=F^0_\om(\phi)$, so we can
apply  Corollary \ref{criterio-del-beta} thus proving the existence of a
\KE metric on $M$.  \fine

The reader will notice that assumption (4) on the ramification
divisors is used only to ensure that \eqref{eq:from-ass-4} holds for
some constant $C_2>0$. This allows to bound
\begin{equation*}
  \intv e^{-(1+\beta)} \om^n \qquad \text{ with } \qquad \intv
e^{-(1+\beta)} \sum_{i=1}^k (\pi_i ^* \om_i)^n . 
\end{equation*}
If the intersection of the ramification divisors is non-vacuous, the
sum of the pull-back measures is degenerate along it.  Nevertheless,
under some numerical assumptions, it is still possible to bound the
integral on the left with the one on the right. 
\begin{prop}\label{criterio-del-c}
  Let $M$ be an $n$-dimensional Fano manifold. Assume that 
  ramified coverings $\pi_i : M \ra M_i$ are given for $i=1, ...,k$,
  satisfying the following assumptions:
   \begin{enumerate}
   \item $M_i$ is a Fano manifold and admits a \KE metric;
   \item  the coverings are Galois, i.e. $M_i=M/G_i$ for some
     finite group $G_i$;
   \item the groups $G_i$ are contained in some \emph{compact}
     subgroup $G \subset \Aut(M)$;
  \item there are (positive) rational numbers $\beta_i$ such that numerically
     \begin{equation*}
       R(\pi_i) = -\beta_i K_M. 
     \end{equation*}
   \end{enumerate}
   Define $\eta\in \cinf(M)$ by
  \begin{equation}
    \label{eq:def-di-eta}
    \frac{1}{k}\sum_{i=1}^k \pi_i^*\om_i^n = \eta \om^n,
  \end{equation}
  and put
  \begin{equation}
    \label{eq:def-di-c}
    c\perdef \sup\{\lambda \geq 0 :  \eta^{-\lamma} \in L^1(M, \om^n)\}
    \end{equation}
and  $\beta \perdef \min \beta_i$. If
\begin{equation}
  \label{ipotesi-su-beta-c}
  \frac{1}{c} < \beta,
\end{equation}
then $M$ admits a \KE metric. 
\end{prop}
\dimo[ of Proposition \ref{criterio-del-c}]
Proceeding as in the proof of Theorem \ref{criterio} one shows that
for any $\phi\in P_G(M,\om)$ 
\begin{equation}
  F_\om^0(\phi) \geq \frac{1}{1+\beta}\log \biggl [ \intv e^{-(1+\beta)\phi} \eta \om^n
  \biggr ] - C_1. 
  \label{pezzo-de-nuevo}
\end{equation}
(Compare with equation \eqref{eq:compare}.) 
It follows from \eqref{ipotesi-su-beta-c} that
we can 
choose a real number $s$ such that
\begin{equation}
  \label{doppio-beta-k}
  1 + \frac{1}{c} < s < 1 +\beta. 
\end{equation}
Put
\begin{equation*}
  \gamma = \frac{1}{s} (1+\beta) - 1. 
\end{equation*}
It follows that $s >1$ and $ \gamma >0$.  Applying H\"{o}lder
inequality with exponent $s$ we see that
\begin{equation}
\begin{gathered}
  \intv e^{-(1+\gamma) \phi} \om^n = \intv e^{-(1+\gamma) \phi}
  \eta^{1/s} \eta^{-1/s} \om^n \leq \\ \leq \biggl [ \intv
  e^{-s(1+\gamma)\phi} \eta \om^n \biggr ]^{\frac{1}{s}} .  \biggl [
  \intv \eta^{-\frac{s'}{s}} \om^n \biggr]^{\frac{1}{s'}}. 
\end{gathered}
\label{holderone} 
\end{equation}
But  \eqref{doppio-beta-k}
\begin{equation*}
  \frac{s'}{s} = \frac{1}{s-1} < c
\end{equation*}
so by the  definition of $c$  
\begin{equation*}
  C_2= \biggl [ \intv \eta^{-\frac{s'}{s}} \om^n \biggr]^{\frac{1}{s'}}
  < +\infty. 
\end{equation*}
On the other hand, $s(1+\gamma)=1+\beta$, so taking the logarithm on both
sides of \eqref{holderone} we get
\begin{equation*}
  \log \biggl [ \intv e^{-(1+\gamma) \phi} \om^n \biggr ] \leq
  \frac{1}{s}\log \biggl [ \intv e^{-(1+\beta)\phi} \eta \om^n \biggr
  ] + \log C_2
\end{equation*}
and applying \eqref{pezzo-de-nuevo}
\begin{equation*}
  F_\om^0(\phi) \geq \frac{s}{1+\beta} \log \biggl [ \intv e^{-(1+\gamma)\phi} \om^n
  \biggr ] - C_3. 
\end{equation*}
Since $\gamma > 0$ we can still apply Corollary
\ref{criterio-del-beta} to get the existence of a \KE metric. 
\fine

It is clear that the last proposition is of some use only if $c$ can be computed or at least bounded from below.              This number is an instance of an interesting invariant of a singularity studied - among others - by Demailly and Koll\'{a}r (see \cite{demailly-kollar-exponent} and \cite{kollar-singularities-pairs}). 
Indeed, 
in the situation of Proposition \ref{criterio-del-c}, let $\mathcal{I}$ be the ideal sheaf on $M$ that on any coordinate chart $U$ is given by $\mathcal{I}=(f_1, ..., f_k)$, where  $f_1, ..., f_k \in \OO_M(U)$ are local defining equations for the divisors $R(\pi_1), ..., R(\pi_k)$.  The \emph{complex singularity exponent} of $\mathcal{I}$ at a point  $x \in U $ is defined as
\begin{equation}
  \label{eq:complex-singularity-exp-def}
  c_x(\mathcal{I}) = \sup \{ \lambda \geq 0 : e^{-2\lambda \phi} \text{ is $L^1$ on a neighbourhood of $x$}\},
\end{equation}
where
\begin{equation}
  \label{eq:def-di-phi-DK}
 \phi = \log (|f_1| + ... + |f_k|) . 
\end{equation}
(See \cite[p. 528]{demailly-kollar-exponent}.) 
Put
\begin{equation}
  \label{eq:sing-exp-globale}
  c_M(\mathcal{I}) = \inf_{x\in M}c_x(\mathcal{I}). 
\end{equation}
\begin{lemma}
  If $c$ is defined by \eqref{eq:def-di-c} and $\mathcal{I}$ is the ideal defined above, then $c=c_M(\mathcal{I})$. 
\end{lemma}
\dimo
Let $(U, z^1 ,...,z^n)$ and $(V, w^1,...,w^n)$ be  coordinate charts on $M$ and  $M_i$ respectively, such that $\pi_i(U)\subset V$. Let $w^s=\pi_i^s(z)$ be the local representation of $\pi_i$. Then the ramification divisor $R(\pi_i)$  is defined by $f_i = \det (\partial \pi_i^s / \partial z^t)$. On the other hand let 
\begin{equation*}
  \begin{split}
    \om & = \I g_{st}\di z^s \wedge \di \barz^t\\
    \om_i & = \I h_{st}\di w^s \wedge \di \bar{w}^t
  \end{split}
\end{equation*}
be the local representations of $\om $ and $\om_i$ on $M$ and $M_i$ respectively. 
It is easy to check that $\pi_i^* \om_i^n = |f_i|^2 \psi_i \om^n$, where
\begin{equation*}
  \psi_i = \frac{| \det(h_{st})|^2}{|\det (g_{st})|^2} . 
\end{equation*}
This is  a smooth positive function, and by restricting $U$ we can assume that 
it  be bounded and uniformly bounded away from $0$.  Cover $M$ with a finite collection of open sets $U_\alpha$ such that this holds for all coverings $\pi_1$, ..., $\pi_k$. On each such $U_\alpha$ we have
$$
  \eta = \frac{1}{k} \bigl ( |f_1|^2 \psi_1 + ... + |f_k|^2 \psi_k\bigr ),
$$
so for some $C>0$
\begin{equation*}
  \frac{1}{C} \eta \leq |f_1|^2 + ... + |f_k|^2 \leq C \eta. 
\end{equation*}
Since $|f_1 | + ... + |f_k| \leq \sqrt{k} \sqrt{ |f_1|^2 + ... + |f_k|^2 } \leq \sqrt{k} \bigl ( |f_1 | + ... + |f_k|\bigr )$, 
the local integrability of $\eta^{-\lambda}$ is equivalent to the local integrability of $e^{-2\lambda \phi}$ (where $\phi$ is defined by \eqref{eq:def-di-phi-DK}).  Taking the minimum over $\alpha$ we get the result. \fine

The complex singularity exponent is in general quite difficult to compute, even for reasonably simple singularities (see \cite[\S 8]{kollar-singularities-pairs}). We present below two cases in which the computation is very simple. Although in many other explicit  examples it is possible to compute   $c$ and to successfully apply Proposition \ref{criterio-del-c}, a general computation of $c$ seems to be hard, although the singularities of the ramification divisors are relatively mild compared to other kinds of singularities.

 We first recall some results on the ramification divisor of a Galois covering. 
\begin{lemma}[Cartan,{ \cite[p. 97]{cartan-pro-lefschetz-suo}}]
\label{cartan-pro-lefschetz-suo}
Given a finite group $G$ acting holomorphically on a complex manifold
$M$ and leaving a point $x\in M$ fixed, there is a biholomorphism
between a neighbourhood of $x$ and a neighbourhood of the origin in
$T_xM$, that intertwines the action of $G$ and the tangent
representation. 
\end{lemma}
\begin{definiz}\label{def-reflectio}
  A \emph{(pseudo)reflection} is a linear map $g \in \Gl(n,\C)$ that
  is diagonalisable and has exactly $n-1$ eigenvalues equal to 1. A
  \emph{reflection group} is a finite subgroup $ G\subset \Gl(n,\C)$
  that is generated by reflections. 
\end{definiz}
The eigenvalues of a reflection $g$ of finite order (i.e. such that
$g^m=1$) are an $m$-th root of unity (with multiplicity 1) and 1 (with
multiplicity $n-1$). (When $m=2$, $g$ is indeed the reflection across its 1-eigenspace.) 
\begin{teo}[Chevalley-Shephard-Todd]
  \label{cavallo-pastore-rospo}
A  finite subgroup   $G \subset \Gl(n, \C)$ is a reflection group if and only if  the affine variety $\C^n/G $ is smooth. 
\end{teo} 
For the proof of this Theorem we refer to \cite[p. 
76]{springer-LNM-585}.  
Let now $\pi:M\ra N=M/G$ be a Galois covering and $x$ a point in $ M$. Denote by $G_x$ the stabiliser. Since the action is properly discontinuous, we can find a neighbourhood $U_x$ of $x$ that is $G_x$-stable and such that $g U_x \cap U_x = \vacuo$ if $g\not \in G_x$. 
 By Cartan's lemma we can assume that $U_x$ be isomorphic to some neighbourhood $V$ of the origin in $T_x M$ with the tangent representation. 
But $U_x/G_x$ is isomorphic to a neighbourhood of $\pi(x)$ in $N$, and therefore is smooth. Hence, Chevalley-Shephard-Todd's theorem implies that $G_x$ acts on $T_x M$ as a reflection subgroup. 
Moreover the invariant theory of finite groups provides a nice model for the map $U_x \ra \pi(U_x)$, and in particular ensures that the projection $\pi$ can be written locally using invariant polynomials: $\pi(z)=(F_1(z), ..., F_n(z))$. Here $F_j$ is a $G_x$-invariant polynomial on $T_x M \cong \C^n$ of degree $d_i$. The polynomial $f=\det (\partial F_i / \partial z^j)$ is a local defining equation for $R(\pi)$. It has degree $(d_1 - 1) + ... + (d_n -1)$. On the other hand the (local) degree of the covering $\pi$ is of course $d_1...d_n$. The inequality $d_1 + ... + d_n -n \leq d_1 ... d_n -1$ implies that in these coordinates the ramification divisor is given by a homogeneous polynomial $f$ whose degree is strictly smaller than the local degree of $\pi$, hence a fortiori smaller than the global degree of $\pi$. 
Thus we have proved the following. 
\begin{lemma}\label{lemmetto-sega-numero-1}
If $\pi : M\ra N$ is a Galois covering between smooth complex manifolds with structure group $G$, then in appropriate coordinate charts around an arbitrary point the local defining equation of  the ramification divisor is a \emph{homogeneous} polynomial of degree less than $\# G$. 
\end{lemma}
The description of the ramification divisor can be made more precise
(see \cite[Exercise 4.3.5 p. 85]{springer-LNM-585}). Let $H$ be a
hyperplane in $\C$. The reflections in $G_x$ that fix $H$ form a
cyclic group. Denote by $e(H)$ its order, and denote by $\ell_H$ a
linear function on $\C^n$ such that $H=\{\ell_H =0\}$. Since there are
a finite number of reflections there are a finite number of
hyperplanes, say $H_1, ..., H_N$, that are fixed by some reflection in
$G_x$. Then on $U_x$ the ramification divisor has the following local
defining equation:
\begin{equation}
  \label{eq:ramificazione-small}
  f=\prod_{i=1}^N \ell_{H_i}^{ e(H_i)-1}=0. 
\end{equation}
If the (reduced) ramification is smooth there is only one hyperplane. 
Since $e(H)\leq \#G_x$, we have proved the following. 
\begin{lemma}\label{lemmetto-sega-numero-2}
  Let $\pi : M\ra N$ be a Galois covering between smooth complex manifolds with structure group $G$. If the reduced divisor associated to the ramification divisor is smooth at $x\in M$, then  there is a holomorphic function $\ell$ defined on some neighbourhood  of $x$, such that $d\ell (x) \neq 0$, and 
$R(\pi)=\{\ell^m=0\}$, with $m\leq \#G - 1$. 
\end{lemma}
We can now give two simple applications of Proposition \ref{criterio-del-c}. 
\begin{teo}\label{un-rivestimento-solo}
  Let $M$ be an $n$-dimensional Fano manifold, and let $\pi :M\ra N$ be a Galois covering with group $G$ onto a \KE manifold $N$. Assume that homologically $       R(\pi) = -\beta K_M$, and that 
\begin{equation}
  d-1 < \beta
  \label{criterio-con-un-beta}
\end{equation}
where $d=\#G = \operatorname{deg}(\pi)$. 
   Then $M$ has a \KE metric. 
\end{teo}
\dimo
Take $x\in M$. If $x$ does not lie in the support of $R=R(\pi)$ then $\eta^{-\lambda}$ is clearly $L^1_{loc}$ for any positive $\lambda$. If $x$ lies in the support of $R$, 
Lemma \ref{lemmetto-sega-numero-1} implies that in appropriate coordinates centered at $x$ the divisor $R$ has a local defining equation that is a homogeneous polynomial $f$ of degree $m$, with $m\leq d-1$. In particular $\operatorname{ord}_x f = m$. The following  general result gives a lower bound for $c_x(f)$ (for the proof see  \cite[Lemma 8.2, p. 438]{open-demailly}). 
\begin{lemma}
  Let $f$ be a holomorphic function on an open set $U\subset \C^n$. If $x\in U$,
then $c_x(f) \geq 1 / \operatorname{ord}_x(f)$. 
\end{lemma}
From this it follows that $c_x(\mathcal{I})=c_x(f) \geq 1/m \geq 1/(d-1)$  for any point of $M$.  Hence $c=c_M(\mathcal{I}) \geq 1/(d-1)$, and an application of Proposition \ref{criterio-del-c} concludes the proof. 
\fine

\begin{teo}\label{criterio2}
  Let $M$ be an $n$-dimensional Fano manifold. Assume that 
  ramified coverings $\pi_i : M \ra M_i$ are given for $i=1, ...,k$,
  satisfying the following assumptions:
   \begin{enumerate}
   \item $M_i$ is a Fano manifold and admits a \KE metric;
   \item the coverings are Galois, i.e. $M_i=M/G_i$;
   \item  the groups $G_i$ are all contained in some fixed \emph{compact}
     subgroup $G \subset \Aut(M)$;
   \item if $V_i $ denotes the reduced divisor of $M$ associated to
     the ramification divisor of $\pi_i$, then the $V_i$'s are smooth
     hypersurfaces, 
     that intersect transversally in a smooth submanifold $V$;
   \item there are (positive) rational numbers $\beta_i$ such that
     \begin{equation*}
       R(\pi_i) = -\beta_i K_M,
     \end{equation*}
     and they satisfy
\begin{equation}
  \frac{1}{d_1 -1 } + ... + \frac{1}{d_k -1 }> \frac{1}{\beta}
\label{ipotesi-sul-beta-k}
\end{equation}
where  $\beta \perdef \min \beta_i$ and $d_i = \# G_i$. 
   \end{enumerate}
   Then $M$ has a \KE metric. 
\end{teo}
\dimo
In order to apply Proposition \ref{criterio-del-c} it is necessary to show that 
\begin{equation}
  \label{eq:a-demontrer-criterio-del-beta-k}
  c \geq \frac{1}{d_1-1} + ... + \frac{1}{d_k-1}. 
\end{equation}
By definition $V=V_1 \cap ... \cap V_k=V(\mathcal{I})_{red}=V(\sqrt{\mathcal{I}})$. 
Let $x$ be a point in $M$. If $x \not\in V$ then $e^{-2\phi(x)} $  is finite (see \eqref{eq:def-di-phi-DK}), and clearly $c_x(\mathcal{I}) = + \infty$. Let $x\in V$.   Using Lemma \ref{lemmetto-sega-numero-2}
we find a neighbourhood $U$ of $x$ and holomorphic functions $\ell_1, ..., \ell_k$ such that $R(\pi_i)=\{\ell_i^{m_i}=0\}$. Since the $V_i$'s cross normally, $\di \ell_1, ... , \di \ell_k$ are linearly independent, hence we can find a coordinate system on a neighbourhood $U$ of $x$ such that $\ell_i=z_i$ for $1\leq i \leq k$. 
Since $m_i\leq d_i -1$, 
in order to prove \eqref{eq:a-demontrer-criterio-del-beta-k} it is enough to show that the integral 
\begin{equation*}
  I(\lamma)=\int_U \bigl ( |z_1^{m_1}| + ... + |z_k^{m_k}| \bigr )^{-2\lambda} \om^n
\end{equation*}
converges for any positive $\lamma < 1/m_1 + ... + 1/m_k$, i.e. that $c_x(\mathcal{I}) \geq 1/m_1 + ... + 1/m_k$. 
 Assuming that the coordinate chart maps $U$ into a polydisk $\Delta^n$ (where $\Delta=\{z\in \C : |z|<1\})$, we get the estimate
\begin{equation*}
\begin{aligned}
  I(\lamma)& \leq C_1\int_{\Delta^n} \frac{1}{\bigl (|z_1|^{m_1} + ... + |z_k|^{m_k}\bigr )^{2\lamma}} \di \mathcal{L}^n =\\
&= C_2 \int_{\Delta^k} \frac{1}{\bigl ( |z_1|^{m_1} + ... + |z_k|^{m_k}\bigr ) ^{2\lamma}} \di \mathcal{L}^k
\end{aligned}
\end{equation*}
$\mathcal{L}^n$ being $2n$-dimensional Lebesgue measure. 
Using polar coordinates in each disk $\Delta$ with $t_i=|z_i|$, we get
\begin{equation*}
  I(\lambda) \leq C_3 \int_0^1 \di t_1 ...\int_0^1 \di t_k \frac{t_1 ... t_k} {(t_1^{m_1} + ... + t_k ^{m_k} )^{2\lamma} }. 
\end{equation*}
With the substitution $t_i = s_i^{1/m_i}$ 
\begin{equation*}
  I(\lamma) \leq C_4 \int_0^{1} \di s_1 ...\int_0^1 \di s_k
\frac{s_1^{\frac{2}{m_1} -1} ... s_k^{\frac{2}{m_k} -1 } }
{(s_1 + ... +s_k)^{2\lamma}}. 
\end{equation*}
If $\lambda < 1/m_1 + ... + 1/m_k$, we can choose $\lamma_1, ..., \lamma_k$ such that $0< \lamma_i < 1/{m_i}$ and $\lamma = \lamma_1 + ... + \lamma_k$. Since $s_1 + ... +s_k \geq s_i$ we get
\begin{equation*}
  I(\lamma) \leq C_4  \int_0^{1} \di s_1 ... \int_0^1 \di s_k
\prod_{i=1}^k
\Biggl [ \frac{s_i^{\frac{2}{m_i} -1}} {(s_1 + ...+s_k)^{2\lamma_i}} 
\Biggr ]
\leq
C_4 \prod_{i=1}^k \int_0^{1} 
s_i ^{2 ( \frac{1}{m_i} -\lamma_i) -1} 
\di s_i . 
\end{equation*}
And this converges since $\frac{1}{m_i} -\lamma_i > 0$ for every $i$. 
\fine

\section{Examples}
\label{esempiozzi-lisci}

Consider the hypersurface
\begin{equation*}
  M=\{ x_0^d + ... +x_{k-1}^d + f(x_k , ..., x_{n+1}) = 0 \} \subset
  \PP^{n+1}
\end{equation*}
where $f$ is any homogeneous polynomial of degree $d$ such that $M$ is
smooth. Note that this is equivalent to saying that
\begin{equation*}
  V = M \cap \{ x_0=...=x_{k-1}=0\} \cong \{f=0\} \subset \PP^{n+1-k}
\end{equation*}
be smooth. 
\begin{prop}\label{ipersupKE}
  If $k > n+2 -d$ then $M$ admits a \KE metric. 
\end{prop}
\dimo $M$ admits $k$ Galois $\ciclico{d}$-coverings $\pi_i : M \ra
\PP^n$ obtained by deleting the $i$-th coordinate, $\pi(x_0, ...,
x_{n+1}) = (x_0, ..., \widehat{x_i}, ..., x_{n+1})$. $G_i=
\ciclico{d}$ acts by multiplication by roots of unity on the $i$-th
coordinate of $\PP^{n+1}$.  $R(\pi_i) = \{
x_i^{d-1}=0\}=\OO(d-1)=-\beta K_M$ with
 \begin{equation*}
  \beta = \frac{d-1}{n+2 -d}. 
\end{equation*}
Since the groups $G_i$ commute, they generate a subgroup of $\Aut(M)$
which is isomorphic to $G_0 \times ... \times G_{k-1}$. Therefore they
all lie inside this (finite) compact subgroup of $\Aut(M)$. The
ramifications are smooth hyperplane sections, and their intersection
is the submanifold $V$ above.  Therefore a straightforward application
of Theorem \ref{criterio2} yields the existence of the \KE metric. 
\fine
\begin{prop}
  Let $M \subset \PP^{n+m}$ be a complete intersection of $m$
  hypersurfaces of degree $d$, given by equations of the form
   \begin{multline*}
    F_j (x_0, ..., x_{n+m}) = a^j_0 x_0^d + ... a^j_{k-1} x_{k-1}^d +
    f_j(x_k,..., x_{n+m}) =0, \\j=1,...,m. 
  \end{multline*}
   I.e. the equations are diagonal in the first $k$ coordinates.  If
  $n+2 -d < k$, then $M$ admits a \KE metric. 
\end{prop}
\dimo We proceed by induction over $m$. For $m=1$ it is the last
Proposition. Let $m>1$, and assume that the result is true for
intersections of $m-1$ hypersurfaces. If we delete one of the first
$k$ coordinates, for example $x_0$, we get a degree $d$ covering
$\pi_0 : M \ra M_0 \subset \PP^{n+m-1}$ over a manifold with equations
\begin{equation*}
\begin{gathered}  a^1_0 F_j - a^j_0 F_1 =
  b^j_1 x_1^d + ... + b^j_{k-1} x_{k-1}^d + h_j(x_k, ..., x_{k+m}) =0
  ,\\ \text{where} \qquad \begin{cases} b^j_s= a^1_0 a^j_s - a^j_0
    a^1_s\\ h_j=a^1_0 f_j - a^j_0 f_1
\end{cases}
\end{gathered}
\end{equation*}
Therefore the base of the covering has equations of the same form, but
in smaller number.  By induction it has a \KE metric. Moreover we can
do the same with any other coordinate $x_1, ..., x_{k-1}$, so we get
$k$ coverings over \KE manifolds. The ramifications are smooth, as
well as their intersection, and
\begin{equation*}
  \beta= \frac{d-1}{n+m+1-md}. 
\end{equation*}
Since $n+1+m(1-d) \leq n+2 -d <k$, we see that $\beta > (d-1) /k$, and
we can apply Theorem \ref{criterio2} to get the \KE metric. \fine

When $d=2$, i.e. when we are intersecting quadrics, one needs $k=n+1$,
which means that all the quadrics are in diagonal form. If $m=2$, the
following result says that in this way we get \emph{ all} the
intersection of two quadrics. 

\begin{teo}\label{segre}
  If $Q_1, Q_2$ are quadrics in $ \PP^{n+2}$, such that their
  intersection $M=Q_1 \cap Q_2$ is smooth and $n$-dimensional, then
  there is a system of homogeneous coordinates $(x_0:...:x_{n+2})$
  such that
  \begin{equation}
    \label{diagonalizzati-da-bravina}
    \begin{split}
      Q_1 &= \{ x_0^2 + ... +x_{n+2}^2 =0\} \\ Q_2 &=\{ \lambda_0
      x_0^2 + ... + \lambda_{n+2}x_{n+2}^2 =0\}
    \end{split}
  \end{equation}
  with $\lambda _i \neq \lambda_j$ for $i\neq j$. 
\end{teo}
For this classical result we refer the reader to the detailed proof
given by Miles Reid in his PhD thesis \cite[p. 36]{reid-thesis}. 

\begin{cor}
  Any smooth intersection of two quadrics $M=Q_1 \cap Q_2$ in
  $\mathbb{P}^{n+2}$ has a \KE metric. 
\end{cor}
Note that this gives the whole moduli space of such manifolds. In fact
a result of Fujita says that these manifolds are characterised by
simple numerical invariants (see \cite[p. 54, Theorem 3.2.5
(iv)]{fano-varieties}).

Browsing through the list of Fano 3-folds with $\rho=h^{1,1}=1$ (see
e.g. \cite[p. 215]{fano-varieties}) we see that some of them are
already defined as coverings. These are the manifolds that Iskovskikh
called \emph{hyperelliptic} because the anticanonical linear system
$|-K_M|$ determines a morphism that is a double cover onto its image
$M'$.  The branching divisor $B\subset M'$ is smooth, and the pairs
$(M',B)$ can be classified.  The possibilities are the following ones
(see \cite[p. 33-34]{fano-varieties}):
\begin{enumerate}\label{enumerazioni}
\item [a)] $M'=\PP^3$ and $B$ is a sextic surface;
\item [b)] $M'=Q^3$ is the 3-dimensional quadric, and $B$ is cut out
  by a quartic surface;
\item [c)] $M' \subset \PP^6$ is a cone over the Veronese surface, and
  $B$ is cut out by a cubic hypersurface. 
\end{enumerate}
Using Theorem \ref{criterio2} we will show that the manifolds in (a)
and (b) admit a \KE metric. Actually the same holds for analogous
coverings in arbitrary dimension.  Whether (c) can be dealt with these
methods is not clear.

\begin{teo}
  Let $M$ be an $n$-dimensional Fano manifold that admits a double
  covering $\pi$ over $\PP^n$ with branching divisor a smooth
  hypersurface of degree $2d$, with $\frac{n+1}{2} < d\leq n$.  Then
  $M$ admits a \KE metric. 
\end{teo}
\dimo That these coverings are smooth depends on the fact that the
branching divisor is smooth, see \cite[p. 42]{barth-peters-vandeven}. 
Recall that given for a double cover  $\pi: M \ra N$, the
ramification and branching divisors are related by
\begin{equation}
  R = \frac{1}{2} \pi^* B. 
\end{equation}
Since $B=\OO(2d)$, $R=\pi^* \OO(d)$ and it follows from Hurwitz
formula that $-K_M = \pi^* \OO(n+1 -d)$. Therefore $R=-\beta
K_M$ with
\begin{equation}
  \beta = \frac{d}{n+1-d}. 
\end{equation}
As the pull-back of an ample line bundle by a finite map is ample,
these manifolds are Fano for $1 \leq d \leq n$.  In order to apply
Theorem \ref{un-rivestimento-solo} we only need to check that
\eqref{criterio-con-un-beta} holds, with $d=2$, i.e.  $1 < \beta$. And
this holds if and only if $d > (n+1)/2$.  \fine
\begin{teo}
  Let $M$ be an $n$-dimensional Fano manifold that is a double cover
  of the quadric $Q_n\subset \PP^n$ ramified along a smooth divisor
  cut out by a hypersurface of degree $2d$, with $\frac{n}{2} < d <
  n$.  Then $M$ admits a \KE metric. 
\end{teo}
\dimo Denote by $\pi : M \ra Q_n$ the covering and by $ i : Q_n
\hookrightarrow \PP^{n+1}$ the inclusion.  Put $\phi = i \pi$. Then
$B=i^*\OO(2d)$, and $R=\demi \pi^*B = \phi^* \OO(d)$, while $-K_M
= \phi^* \OO(n) - \phi^*\OO(d)=\phi^*\OO(n-d)$.  So for $d<n$ $M$ is
Fano and
\begin{equation*}
  \beta= \frac{d}{n-d}. 
\end{equation*}
When $2d >n$, $\beta >1 $, and Theorem \ref{un-rivestimento-solo}
yields the result.  \fine In case $n=3$, $d$ has to be equal to $2$,
i.e. $B$ is cut out by a quartic, and the branching divisor is an
octic hypersurface contained in $Q_3$. When $d=1$ (and $n$ arbitrary),
it is not possible to apply Theorem \ref{criterio2}, but in this case
the manifold $M$ is simply the intersection of two quadrics.

\end{document}